\documentclass[11pt]{article}
\usepackage{geometry} 
\usepackage{latexsym,amsmath,stackrel,color,amsthm,amssymb,epsfig,graphicx,mathrsfs}
\usepackage{graphicx}
\usepackage{amssymb}
\usepackage[linktocpage=true]{hyperref}
\usepackage{setspace}
\usepackage{multirow}
\usepackage{amssymb, amsmath, amsthm, graphicx,mathrsfs}

\usepackage{indentfirst,subfig}
\usepackage{lipsum}

\usepackage{fancyvrb}
\usepackage{verbatim}

\usepackage{moreverb}

\usepackage{graphicx}
\usepackage{amssymb}
\usepackage{stmaryrd}
\usepackage{enumitem}
\usepackage{mathtools}

\usepackage{amsmath}
\usepackage{systeme}
\usepackage[T1]{fontenc}
\usepackage{babel}
\usepackage{caption}
\usepackage{subcaption}
\usepackage[dvipsnames]{xcolor}
\usepackage{tkz-graph}

\newcommand{\pushright}[1]{\ifmeasuring@#1\else\omit\hfill$\displaystyle#1$\fi\ignorespaces}

\usepackage[numbers,sort&compress]{natbib}
\usepackage{indentfirst,graphics,epsfig}
\usepackage{graphicx}
\usepackage{graphics}
\usepackage{graphicx,psfrag}
\usepackage{graphics}
\usepackage{tikz}
\usetikzlibrary{trees,arrows}
\usetikzlibrary{automata,arrows,positioning,calc}
\usepackage {xcolor}
\definecolor {processblue}{cmyk}{0.10,5,0,0}
\usepackage{pgf,tikz}
\usepackage{mathrsfs}
\usetikzlibrary{arrows}

\usepackage{latexsym,amsmath,amssymb,amsfonts,epsfig,psfrag,url,graphics,ifpdf,multicol}
\usepackage{color,colordvi,amscd,amsthm} 
\usepackage{tikz}
\usepackage[numbers,sort&compress]{natbib}
\usepackage{indentfirst,graphics,epsfig}
\usepackage{graphicx}
\usepackage{graphics}
\usepackage{graphicx,psfrag}
\usepackage{graphics}

\usepackage{tabularx}
\usepackage{booktabs}
\usepackage{floatrow}
\usepackage{multirow}
\usepackage{graphicx}
\usepackage{caption}
\usepackage[ruled]{algorithm2e}
\newfloatcommand{capbtabbox}{table}[][\FBwidth]
\usepackage{float}

\usepackage{indentfirst,subfig}

\newtheorem{theorem}{Theorem}[section]
\newtheorem{lemma}[theorem]{Lemma}

\newtheorem{example}[theorem]{Example}

\newtheorem{hypothesize}[theorem]{Hypothesize}

\usepackage{blindtext}
\usepackage{titlesec}
\usepackage{tikz}
\usetikzlibrary {decorations.pathmorphing,shadows}
\usepackage{neuralnetwork} 
\usepackage{pgf,tikz}
\usepackage{mathrsfs}
\usepackage{bm}
\usetikzlibrary{arrows, graphs,graphs.standard, positioning}
\usepackage{blindtext}
\usepackage{adjustbox}
\usetikzlibrary {graphs}
\usetikzlibrary {graphs.standard}
\usepackage[all]{xy}
\usepackage{lipsum} 
\usepackage[]{mdframed}

\usepackage{array}
\newcolumntype{P}[1]{>{\centering\arraybackslash}p{#1}}
\newcolumntype{M}[1]{>{\centering\arraybackslash}m{#1}}
\begin{document}
\begin{flushleft}
UDC: 519.172.1  \hspace{8cm} $66^Th$ All-Russia Scientific Conference MIPT.
\end{flushleft}
\begin{center}
\textbf{  }\\
	\textbf{Sigma index in Trees with Given Degree Sequences}
\end{center}
\begin{center}
	Jasem Hamoud$^1$, Artem Kurnosov$^2$ \\
	$^1,^2$ Moscow Institute of Physics and Technology 
	(National Research University).
	
\end{center}
\section{Introduction.}

   The sigma index in graph theory refers to a measure of the degree differences between vertices in a graph. The goal is to determine the graphs that have the maximum sigma index within certain classes of graphs. Abdo, Dimitrov, and Gutman characterized the graphs with the greatest sigma index among all connected graphs of a fixed order in\cite{1} as: 
   \begin{enumerate}
   	\item  	$\sigma \left( G \right) = irr\left( G \right) = \sum\limits_{uv \in E(G)}^{} {\left| {{d_G}(u) - {d_G}(v)} \right|}  = 0$, if and only if G is regular, then that’s mean: ${d_G}(u) = {d_G}(v);\forall uv \in E(G)$;
   	\item  	$\sigma \left( G \right) = irr\left( G \right) = \sum\limits_{uv \in E(G)}^{} {\left| {{d_G}(u) - {d_G}(v)} \right|}  > 0$, if and only if G is stepwise regular, then that’s mean: ${d_G}(u) > {d_G}(v);\forall uv \in E(G){\rm{ or }}{d_G}(v) > {d_G}(u);\forall uv \in E(G)$.
   \end{enumerate}
Also, Gutman, Togan, Yurttas, Cevik, Cangul introduced sigma index $\sigma$  in graph in \cite{2} as: 

\[\sigma  = \sum\limits_{uv \in E\left( G \right)}^{} {{{\left( {{d_G}\left( u \right) - {d_G}\left( v \right)} \right)}^2}} .\]

\section{Main Result.}

 \begin{hypothesize}
	let be a sequence $\tilde{d}=(d_1,d_2,d_3); d_1\geq d_2\geq d_3$ of order $n$ where $n\leq 4$ we have: $irr_{max}-irr_{min}=2(d_2-d_3)$ that is mean: $d(irr_{max},irr_{min})=2(d_2-d_3)$ if $d_2=d_3 \Longrightarrow d_G(irr_{max},irr_{min})=2$, where The distance $d_G(u,v)$ between two vertices $u$, $v$ of $G$ is the length of one of the shortest $(u,v)$-path in $G$
\end{hypothesize}
\begin{lemma}
	A sequence of positive integers ${d_1},{d_2},{d_3};{\rm{ }}{d_1} \ge {d_2} \ge {d_3}$  and defined numbers : $a\geq b\geq c$. so that, sigma index defined as: 
	\[
	irr=\left\lbrace
	\begin{aligned}
		& irr_{max}=(a-1)^2+(b-1)^2 +(c-1)(c-2)(a-c)(b-c) \\
		& irr_{min}=(a-1)^2+(c-1)^2+(b-1)(b-2)+(a-c)
	\end{aligned} 
	\right. 
	\]
	we have 3 cases as: 
	\begin{table}[H]
		\centering
		\begin{tabular}{ |p{4cm}|p{4cm}|p{4cm}| }
			\hline
			\hspace{1cm} case 1 & \hspace{1cm} case 2 & \hspace{1cm} case 3 \\
			\hline
			\begin{tikzpicture}[scale=1]
				\draw [line width=2pt] (3,4)-- (2,2);
				\draw [line width=2pt] (3,4)-- (4,2);
				\draw (3.025660215512478,4.746014977092217) node[anchor=north west] {$d_1$};
				\draw (1.3851544958474011,2.486450495289372) node[anchor=north west] {$d_3$};
				\draw (4.04710717303979,2.1923975832739333) node[anchor=north west] {$d_2$};
				\begin{scriptsize}
					\draw [fill=black] (3,4) circle (2.5pt);
					\draw [fill=black] (2,2) circle (2.5pt);
					\draw [fill=black] (4,2) circle (2.5pt);
				\end{scriptsize}
			\end{tikzpicture} 			&  			\begin{tikzpicture}[scale=1]
				
				\draw [line width=2pt] (7,4)-- (6,2);
				\draw [line width=2pt] (7,4)-- (8,2);
				\draw (6.973893188674715,4.65362045684315) node[anchor=north west] {$d_2$};
				\draw (5.697084491765575,2.065954831107291) node[anchor=north west] {$d_3$};
				\draw (8.046412494078393,2.0829789470660796) node[anchor=north west] {$d_1$};
				\begin{scriptsize}
					\draw [fill=black] (7,4) circle (2.5pt);
					\draw [fill=black] (6,2) circle (2.5pt);
					\draw [fill=black] (8,2) circle (2.5pt);
				\end{scriptsize}
			\end{tikzpicture} &
			\begin{tikzpicture} [scale=1]
				\draw [line width=2pt] (11,4)-- (10,2);
				\draw [line width=2pt] (11,4)-- (12,2);
				\draw (10.906463975154866,4.7727892685546704) node[anchor=north west] {$d_3$};
				\draw (9.970137597421497,1.9638101353545598) node[anchor=north west] {$d_1$};
				\draw (11.944935048640968,1.9127377874781941) node[anchor=north west] {$d_2$};
				\begin{scriptsize}
					\draw [fill=black] (11,4) circle (2.5pt);
					\draw [fill=black] (10,2) circle (2.5pt);
					\draw [fill=black] (12,2) circle (2.5pt);
				\end{scriptsize}
			\end{tikzpicture}\\
			\hline
		\end{tabular}
		\caption{positive integer ${d_1},{d_2},{d_3}$.}
	\end{table}

	we have for Albertson Index irregularity : $case 1 = case 2 \leq case 3$.
\end{lemma}

\begin{example}
	let a sequence ${d_1},{d_2},{d_3}{\rm{   ; }}{d_1} \ge {d_2} \ge {d_3}$ and we will suppose sequences is: $(4,3,2)$ we have 3  cases as: 
	
	\begin{table}[H]
		\centering
		\begin{tabular}{ |p{6cm}|p{6cm}|p{6cm}| }
			\hline
			\hspace{2cm} case 1 & \hspace{2cm} case 2 & \hspace{2cm} case 3 \\
			\hline
			\begin{tikzpicture}[scale=1]
				\draw [line width=2pt] (4,6)-- (2,4);
				\draw [line width=2pt] (4,6)-- (4,4);
				\draw [line width=2pt] (4,6)-- (6,4);
				\draw [line width=2pt] (2,4)-- (0.6133153927217206,2.6467527994704714);
				\draw [line width=2pt] (2,4)-- (1.9645890254396188,2.5466584563061825);
				\draw [line width=2pt] (2,4)-- (2.8654381139182177,2.84694148579905);
				\draw [line width=2pt] (4,4)-- (4.016523060307538,2.7468471426347607);
				\draw (3.9164287171432495,7.150998241863484) node[anchor=north west] {$3$};
				\draw (1.5141644812003194,4.92389910645805) node[anchor=north west] {$4$};
				\begin{scriptsize}
					\draw [fill=black] (4,6) circle (2.5pt);
					\draw [fill=black] (2,4) circle (2.5pt);
					\draw [fill=black] (4,4) circle (2.5pt);
					\draw [fill=black] (6,4) circle (2.5pt);
					\draw [fill=black] (0.6133153927217206,2.6467527994704714) circle (2.5pt);
					\draw [fill=black] (1.9645890254396188,2.5466584563061825) circle (2.5pt);
					\draw [fill=black] (2.8654381139182177,2.84694148579905) circle (2.5pt);
					\draw [fill=black] (4.016523060307538,2.7468471426347607) circle (2.5pt);
				\end{scriptsize}
			\end{tikzpicture} 			&  			\begin{tikzpicture}[scale=1]
				
				\draw [line width=2pt] (4,6)-- (2,4);
				\draw [line width=2pt] (4,6)-- (4,4);
				\draw [line width=2pt] (4,6)-- (6,4);
				\draw [line width=2pt] (2,4)-- (0.6133153927217206,2.6467527994704714);
				\draw [line width=2pt] (2,4)-- (1.9645890254396188,2.5466584563061825);
				\draw [line width=2pt] (4,4)-- (4.016523060307538,2.7468471426347607);
				\draw (3.9164287171432495,7.150998241863484) node[anchor=north west] {$4$};
				\draw (1.5141644812003194,4.92389910645805) node[anchor=north west] {$3$};
				\draw [line width=2pt] (4,6)-- (1.5391880669913907,5.849771780727725);
				\begin{scriptsize}
					\draw [fill=black] (4,6) circle (2.5pt);
					\draw [fill=black] (2,4) circle (2.5pt);
					\draw [fill=black] (4,4) circle (2.5pt);
					\draw [fill=black] (6,4) circle (2.5pt);
					\draw [fill=black] (0.6133153927217206,2.6467527994704714) circle (2.5pt);
					\draw [fill=black] (1.9645890254396188,2.5466584563061825) circle (2.5pt);
					\draw [fill=black] (4.016523060307538,2.7468471426347607) circle (2.5pt);
					\draw [fill=black] (1.5391880669913907,5.849771780727725) circle (2.5pt);
				\end{scriptsize}
			\end{tikzpicture} &
			\begin{tikzpicture} [scale=1]
				\draw [line width=2pt] (3.9414523029343216,5.849771780727725)-- (2.915485285500362,4.448450976427677);
				\draw [line width=2pt] (3.9414523029343216,5.849771780727725)-- (5.142584420905787,4.248262290099098);
				\draw [line width=2pt] (2.915485285500362,4.448450976427677)-- (1.714353167528897,2.947035828963339);
				\draw [line width=2pt] (2.915485285500362,4.448450976427677)-- (3.315862658157517,2.7468471426347607);
				\draw [line width=2pt] (5.142584420905787,4.248262290099098)-- (4.567041947711126,2.84694148579905);
				\draw (2.2148248833503406,5.4994415796527125) node[anchor=north west] {$4$};
				\draw (5.517938207771869,5.474417993861641) node[anchor=north west] {$3$};
				\draw [line width=2pt] (5.142584420905787,4.248262290099098)-- (5.86826840884688,2.8969886573811943);
				\draw [line width=2pt] (2.915485285500362,4.448450976427677)-- (2.615202256007495,2.696799971052616);
				\begin{scriptsize}
					\draw [fill=black] (3.9414523029343216,5.849771780727725) circle (2.5pt);
					\draw [fill=black] (2.915485285500362,4.448450976427677) circle (2.5pt);
					\draw [fill=black] (5.142584420905787,4.248262290099098) circle (2.5pt);
					\draw [fill=black] (1.714353167528897,2.947035828963339) circle (2.5pt);
					\draw [fill=black] (3.315862658157517,2.7468471426347607) circle (2.5pt);
					\draw [fill=black] (4.567041947711126,2.84694148579905) circle (2.5pt);
					\draw [fill=black] (5.86826840884688,2.8969886573811943) circle (2.5pt);
					\draw [fill=black] (2.615202256007495,2.696799971052616) circle (2.5pt);
			\end{scriptsize}\end{tikzpicture}\\
			\hline
			\hspace{2cm}	$irr=14$ & \hspace{2cm} $irr=14$ & \hspace{2cm} $irr=16$ \\
			\hline
			\hspace{2cm} $\sigma=34$ & \hspace{2cm} $\sigma=32$ & \hspace{2cm} $\sigma=40$ \\
			\hline
			
		\end{tabular}
		\caption{Albertson Index and Sigma Index for  a sequence ${d_1},{d_2},{d_3}{\rm{   ; }}{d_1} \ge {d_2} \ge {d_3}$.}
	\end{table}	
\end{example}

\begin{example}
	let a sequence ${d_1},{d_2},{d_3},{d_4}{\rm{   ; }}{d_1} \ge {d_2} \ge {d_3} \ge {d_4}$ and we will suppose sequences is: $(8,5,4,2)$ we have 12 options, for example we take: 
	
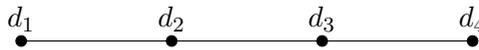
\begin{figure}[H]
		\centering
		\begin{tikzpicture}[scale=1]
			\filldraw 
			(0,0) circle (2pt) node[align=left,   above] {$d_1$} --(2,0) circle (2pt) node[align=right, above] {$d_2$} -- (4,0) circle (2pt) node[align=right, above] {$d_3$}--(6,0) circle (2pt) node[align=right, above] {$d_4$};
		\end{tikzpicture}
	\end{figure}
	Let we take:  
	
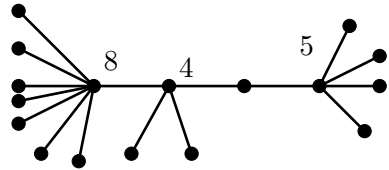
\begin{figure}[H]
		\centering
		\begin{tikzpicture}[scale=1] 
			\draw [line width=1pt] (3,3)-- (2,4);
			\draw [line width=1pt] (3,3)-- (2,3.5);
			\draw [line width=1pt] (3,3)-- (2,3);
			\draw [line width=1pt] (3,3)-- (2,2.8);
			\draw [line width=1pt] (3,3)-- (2,2.5);
			\draw [line width=1pt] (3,3)-- (2.3,2.1);
			\draw [line width=1pt] (3,3)-- (2.8,2);
			\draw [line width=1pt] (3,3)-- (4,3);
			\draw [line width=1pt] (4,3)-- (5,3);
			\draw [line width=1pt] (4,3)-- (3.5,2.1);
			\draw [line width=1pt] (4,3)-- (4.3,2.1);
			\draw [line width=1pt] (5,3)-- (6,3);
			\draw [line width=1pt] (6,3)-- (6.4,3.8);
			\draw [line width=1pt] (6,3)-- (6.8,3.4);
			\draw [line width=1pt] (6,3)-- (6.8,3);
			\draw [line width=1pt] (6,3)-- (6.6,2.4);
			\draw (3,3.6) node[anchor=north west] {$8$};
			\draw (4,3.5) node[anchor=north west] {$4$};
			\draw (5.6,3.8) node[anchor=north west] {$5$};
			\begin{scriptsize}
				\draw [fill=black] (3,3) circle (2.5pt);
				\draw [fill=black] (2,4) circle (2.5pt);
				\draw [fill=black] (2,3.5) circle (2.5pt);
				\draw [fill=black] (2,3) circle (2.5pt);
				\draw [fill=black] (2,2.8) circle (2.5pt);
				\draw [fill=black] (2,2.5) circle (2.5pt);
				\draw [fill=black] (2.3,2.1) circle (2.5pt);
				\draw [fill=black] (2.8,2) circle (2.5pt);
				\draw [fill=black] (4,3) circle (2.5pt);
				\draw [fill=black] (5,3) circle (2.5pt);
				\draw [fill=black] (3.5,2.1) circle (2.5pt);
				\draw [fill=black] (4.3,2.1) circle (2.5pt);
				\draw [fill=black] (6,3) circle (2.5pt);
				\draw [fill=black] (6.4,3.8) circle (2.5pt);
				\draw [fill=black] (6.8,3.4) circle (2.5pt);
				\draw [fill=black] (6.8,3) circle (2.5pt);
				\draw [fill=black] (6.6,2.4) circle (2.5pt);
			\end{scriptsize}
		\end{tikzpicture}
	\end{figure}
	We can discuss many cases according to the following: 
	
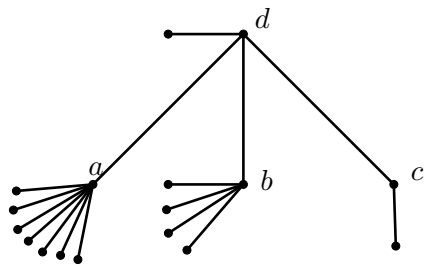
\begin{figure}[H]
		\centering
		\begin{tikzpicture}[scale=1] 
			\draw [line width=1pt] (6,5)-- (4,3);
			\draw [line width=1pt] (6,5)-- (6,3);
			\draw [line width=1pt] (6,5)-- (8,3);
			\draw [line width=1pt] (4,3)-- (2.9818652942587778,2.912736742127259);
			\draw [line width=1pt] (4,3)-- (2.9420336490474974,2.658810003905346);
			\draw [line width=1pt] (4,3)-- (3,2.4);
			\draw [line width=1pt] (4,3)-- (3.1411918751038996,2.245556684838311);
			\draw [line width=1pt] (4,3)-- (3.330392189857482,2.1011669709474194);
			\draw [line width=1pt] (4,3)-- (3.569382061125165,2.056356370084729);
			\draw [line width=1pt] (4,3)-- (3.8,2);
			\draw [line width=1pt] (6,3)-- (5,3);
			\draw [line width=1pt] (6,3)-- (4.978426510474212,2.663788959556756);
			\draw [line width=1pt] (6,3)-- (5.002896977194556,2.3511319237990738);
			\draw [line width=1pt] (6,3)-- (5.249496121010085,2.125976183793591);
			\draw [line width=1pt] (8,3)-- (8.026416914411044,2.179584693318706);
			\draw [line width=1pt] (6,5)-- (5,5);
			\draw (6.016768583055512,5.4879402972858715) node[anchor=north west] {$d$};
			\draw (3.806697611685776,3.4320603239186784) node[anchor=north west] {$a$};
			\draw (6.093864082056783,3.3292663252503187) node[anchor=north west] {$b$};
			\draw (8.089780889534103,3.3635309914731053) node[anchor=north west] {$c$};
			\begin{scriptsize}
				\draw [fill=black] (6,5) circle (1.5pt);
				\draw [fill=black] (4,3) circle (1.5pt);
				\draw [fill=black] (6,3) circle (1.5pt);
				\draw [fill=black] (8,3) circle (1.5pt);
				\draw [fill=black] (2.9818652942587778,2.912736742127259) circle (1.5pt);
				\draw [fill=black] (2.9420336490474974,2.658810003905346) circle (1.5pt);
				\draw [fill=black] (3,2.4) circle (1.5pt);
				\draw [fill=black] (3.1411918751038996,2.245556684838311) circle (1.5pt);
				\draw [fill=black] (3.330392189857482,2.1011669709474194) circle (1.5pt);
				\draw [fill=black] (3.569382061125165,2.056356370084729) circle (1.5pt);
				\draw [fill=black] (3.8,2) circle (1.5pt);
				\draw [fill=black] (5,3) circle (1.5pt);
				\draw [fill=black] (4.978426510474212,2.663788959556756) circle (1.5pt);
				\draw [fill=black] (5.002896977194556,2.3511319237990738) circle (1.5pt);
				\draw [fill=black] (5.249496121010085,2.125976183793591) circle (1.5pt);
				\draw [fill=black] (8.026416914411044,2.179584693318706) circle (1.5pt);
				\draw [fill=black] (5,5) circle (1.5pt);
			\end{scriptsize}
		\end{tikzpicture}
	\end{figure}
	
	if $d=2$ we have:\\ 
	$
	49 + 16 + 4 + 2 + 9 = 80 \Rightarrow {irr} = 80 $ \\
	$
	343+16+27+84=470 \Rightarrow \sigma =470
	$\\
	From this case we can formed many methods to obtain on maximum value, we will take for example $d=4$ we have;
	
	\[\begin{array}{l}
		3 + 1 + 16 + 49 + 2 + 1 + 4 = 76 \Rightarrow irr  = 76 \\
		9+22+343+64 438 \Rightarrow \sigma =438
	\end{array}\]
	
	We can formulate  that as: \\
	$
	irr={\left( {a - 1} \right)^2} + {\left( {b - 1} \right)^2} + {\left( {c - 1} \right)^2} + \left( {d - a} \right) + \left( {d - b} \right) + \left( {d - c} \right) + \left( {d - 1} \right)\left( {d - 3} \right), \\
	if (d=8,a=5,b=4,c=2) \Rightarrow irr= 74
	$
	\\
	\\
	If $d=5$ we have:  \\
	$
	8+10+49+10=74 \Rightarrow irr=74 \\
	32+19+343+28=422 \Rightarrow \sigma =422
	$ 
	\\
	\\
	If $d=8$ we have:  \\
	$
	35+10+29=74 \Rightarrow irr=74 \\
	64+27+36+26+245=398 \Rightarrow \sigma =398
	$ 
	\\
	we will get on:
	\begin{table}[H]
		
		\centering
		\begin{tabular}{ |M{2.5cm}|M{2.5cm}|M{2.5cm}| }
			\hline
			sequence of d &  irr &  $\sigma$ \\
			\hline
			2	&  	80	 & 470   \\
			\hline
			4		&  76		 &   438 \\
			\hline
			5		&  	74	 &  422  \\
			\hline
			
			8		&  	74	 &  398  \\
			\hline
			
		\end{tabular}
	\end{table}	
	
	In the end we get on: \\
	\begin{table}[H]
		
		\centering
		\begin{tabular}{ |M{2.5cm}|M{2.5cm}|M{2.5cm}|M{2.5cm}|  }
			\hline
			\multicolumn{2}{|c|}{$irr$} &  \multicolumn{2}{|c|}{$\sigma$} \\
			\hline
			$\max$ & $\min$ & $\max$ & $\min$ \\
			\hline
			80   & 74    & 470 &   398\\
			\hline
		\end{tabular}
	\end{table}	
		\end{example}


\begin{thebibliography}{10}
\bibitem{1}
	Hosam, A., Darko, D., Ivan G., 2018, Graphs with maximal $\sigma$ irregularity, Elsevier B.V, doi.org/10.1016/j.dam.2018.05.013.
\bibitem{2}
  	I. Gutman, M. Togan, A. Yurttas, A.S. Cevik, I.N. Cangul, Inverse problem for sigma index, MATCH Commun. Math. Comput. Chem. 79 (2018) 491–508. 	
	
\end{thebibliography}
\end{document}